% gjacobimagnus.tex           Feb. 1994
%  this is Plain TeX
%%
%% Here is a character listing to check to be sure that no
%% unwanted translations took place within the bowels of the net.
%% Upper case letters: ABCDEFGHIJKLMNOPQRSTUVWXYZ
%% Lower case letters: abcdefghijklmnopqrstuvwxyz
%% round parentheses, square brackets, curly braces: ()  []  {}
%% Exclaim, at, sharp, dollar, percent: ! @ # $ %
%% Caret, ampersand, star, underscore, hyphen: ^ & * _ -
%% vertical bar, backslash, tilde, backprime, plus: | \ ~ ` +
%% plus, equal, prime, quote, colon: + = ' " :
%% less than, greater than, slash, question, comma: < > / ? ,
%% period, semicolon: . ;
%%
%%

 \magnification=1200
 \pretolerance=999 \tolerance=1999
 \parindent=5pt

\font\bfx=cmbx10 scaled 1095
\font\ttt=cmtt10 scaled 913
\font\rmt=cmr10 scaled 913

\def\dfrac#1#2{{\displaystyle{#1\over#2}}}
\def\rondsur{\mathaccent"7017 }
       % phantom
\def\odd{{\rm \,odd\,}}
\overfullrule=0pt

% From gentle.tex :
%%%%%%%%%%%%%%%%%%%%%%% headline routines %%%%%%%%%%%%%%%%%%%%%%%%%%%%
\def\gentleheadline{%
\centerline{\rm Simplest generalized Jacobi polynomials.}
}

\newif \iftitlepage    \titlepagetrue
\headline=
  {\iftitlepage \hfil \global\titlepagefalse \else \gentleheadline \fi}
%%%%%%%%%%%%%%%%%%%%%%%%%%%%%%%%%%%%%%%%%%%%%%%%%%%%%%%%%%%%%%%%%%%%%

% verbatim macro for plain TeX (TeXbook p. 380-381 )
% skips blank lines and fails with ! ? and tabs
\def\uncatcodespecials{\def\do##1{\catcode`##1=12}\dospecials}
{\obeyspaces\global\let =\ }
\def\setupverbatim{\ttt \obeylines\uncatcodespecials\obeyspaces
                   \parskip=-3pt
                   \everypar{\qquad} }
\def\verblist#1{\par\begingroup\setupverbatim\input#1 \endgroup}

%\nopagenumbers

\centerline{\bfx
 Asymptotics for the simplest generalized Jacobi polynomials recurrence}

\centerline{\bfx
  coefficients from Freud's equations: numerical explorations.}

\bigskip
\centerline{Alphonse P. Magnus}

\centerline{ Institut Math\'ematique,  Universit\'e Catholique de Louvain}

\centerline{Chemin du Cyclotron 2}

\centerline { B-1348 Louvain-la-Neuve}

\centerline{ Belgium}

\centerline{E-mail: {\tt magnus@anma.ucl.ac.be}}

\bigskip

% if this happens to be accepted for the Gatteschi proceedings:

Dedicated to L. Gatteschi on the occasion of his
                                          $70^{\hbox{\rmt th}}$ birthday

\bigskip

{\bf Abstract.\/} Generalized Jacobi polynomials are orthogonal
 polynomials related to a weight function which is smooth and positive
 on the whole interval of orthogonality up to a finite number of points,
  where algebraic singularities occur.
  The influence of these singular points on the asymptotic behaviour of
 the recurrence coefficients is investigated.

\bigskip

{\bf AMS(MOS) subject classification.\/} 42C05.

\bigskip


\bigskip
{\bf 1.Weight singularities and recurrence coefficients.}

\bigskip
The orthonormal polynomials $p_n(x)=\gamma_n x^n+\cdots$ related to
the weight $w$ satisfy the three-terms recurrence relation
$$a_{n+1}p_{n+1}(x)=(x-b_n)p_n(x)-a_np_{n-1}(x), \eqno(1) $$
with $a_0p_{-1}(x)\equiv0$.

Let $w(x)>0$ hold almost everywhere on the support $[-1,1]$, then one
knows (since 1977) that the recurrence coefficients have limits
$a_n\to 1/2$ and $b_n\to 0$ when $n\to\infty$ (see for instance the survey
in [Nev2]).

Features of $w$ can somehow be ``read'' in the sequence of the recurrence
coefficients. This goes back to Stieltjes and is current practice in
solid-state physics [LaG]. See [Ap] for the case of endpoints
singularities.
 J.P.~Gaspard once showed me a paper by C.Hodges
[Ho] describing the influence of a mild
{\it interior\/} singularity of the form
$$\eqalign{
 w(x) &\sim w(x_0)+ A (x-x_0)^\gamma\,,\qquad x\to x_0, x>x_0,  \cr
      &\sim w(x_0)+ B (x_0-x)^\gamma\,,\qquad x\to x_0, x<x_0,  \cr
          } \eqno(2)
$$ with
$x_0\in (-1,1), 0<w(x)<\infty$ on $(-1,1)$ and $\gamma>0$
(Van Hove singularity [Ho,Mart]) as
$$ a_n-1/2 \sim \xi n^{-\gamma-1}\cos(2n\theta_0-\eta),\quad
   b_n     \sim 2\xi n^{-\gamma-1}\cos((2n+1)\theta_0-\eta),\quad
   n\to\infty
$$
with $x_0=\cos\theta_0$.
I gave a lengthy proof of this in [Mag1],
showing how $\xi$ and $\eta$ are related through Toeplitz determinants
to the Szeg\H o function of
$w(cos\theta)\sin\theta$ (the analytic function $D(z)$ in $|z|<1$,
with $D(0)>0$, without zero in $|z|<1$, and such that the boundary
values satisfy $|D(e^{i\theta})|^2= w(\cos\theta)\sin\theta$) by
$$ \xi e^{i\eta} = \dfrac{\left( \dfrac{\sin\theta_0}2\right)^{\gamma+1}
                            \Gamma(\gamma+1)}
                         {2\pi w(x_0)}
               \,e^{i(\theta_0-2\arg(D(e^{i\theta_0}))}
              \left[ A e^{i\pi\gamma/2} - B e^{-i\pi\gamma/2}\right].
$$
Hodges argument is based on the continued
fraction expansion of
$$f(z)=\int_{-1}^1 \dfrac{w(t)\,dt}{z-t} =
    \dfrac{\mu_0}{ z-b_0 -\dfrac{a_1^2}{z-b_1-\cdots}  }$$
  when $z=x+i\varepsilon$
for small $\varepsilon>0$ and $-1<x<1$,
but the fastest explanation is probably related to
inverse scattering techniques  ( [VA] p.117,  [NV] and references
therein):
let us consider the function $q_n$ defined by the integral
$q_n(z)=\int_{-1}^1 (z-t)^{-1}p_n(t)w(t)dt$ when $z\notin [-1,1]$.
Remark that $f(z)=\sqrt{\mu_0} q_0(z)$.
From (1), $a_{n+1}q_{n+1}(z)=(z-b_n)q_n(z)-a_nq_{n-1}(z)
-\sqrt{\mu_0}\delta_{n,0}$, (with $a_0q_{-1}(z)\equiv 0$), which can
be written as
$$ q_{n+1}(z)-2zq_n(z)+q_{n-1}(z)=
   q_{n+1}(z)-\rho(z)^{-1}q_n(z)-\rho(z)[q_n(z)-\rho(z)^{-1}q_{n-1}(z)]=
   \varepsilon_n(z),\, n=0,1,\ldots$$
with $\varepsilon_n(z)=(1-2a_{n+1})q_{n+1}(z)-2b_n q_n(z)
  +(1-2a_n)q_{n-1}(z)-2\sqrt{\mu_0}\delta_{n,0}$, and where
$\rho(z)$ is the determination of $(z+\sqrt{z^2-1})/2$ such that
$|\rho(z)|>1$ when $z\notin [-1,1]$. After a simple summation,
$\rho(z)^{-N}(q_{N+1}(z)-\rho(z)^{-1}q_N(z))-\rho(z)q_0(z)=
  \sum_0^N \rho(z)^{-n}\varepsilon_n(z)$. While $z\notin [-1,1]$ and
when $N\to\infty$, $\rho(z)^{-N}$ and $q_N(z)\to 0$, so
$$q_0(z)=\dfrac{f(z)}{\sqrt{\mu_0}} =
    -\sum_0^\infty \dfrac{\varepsilon_n(z)}{\rho(z)^{n+1}} =
 \dfrac{2\sqrt{\mu_0}}{\rho(z)} +
    \sum_0^\infty \dfrac{(2a_n-1)(q_{n-1}(z)+\rho(z)q_n(z))
                         +2b_n q_n(z) }{\rho(z)^{n+1}}$$
From the known asymptotic behaviour
$q_n(z)\sim \dfrac{2(2\pi)^{1/2}D(\rho(z)^{-1})}{(\rho(z)-\rho(z)^{-1})
   \rho(z)^n}$ for large $n$ [Ba], we have a first approximate expansion
of $(2\pi)^{-1/2}D(\rho(z)^{-1})(\rho(z)-\rho(z)^{-1})q_0(z)$
 as a series of negative powers of $\rho(z)$:
 $$(2\pi)^{-1/2}D(\rho(z)^{-1})(\rho(z)-\rho(z)^{-1})q_0(z) \approx
    \sum_0^\infty \dfrac{2(2a_n-1)}{\rho(z)^{2n}}+
                  \dfrac{ 2b_n}{\rho(z)^{2n+1}}$$
Under sufficiently strong conditions , this remains valid when
$z=x+i\varepsilon, \varepsilon\to 0, \varepsilon>0,
-1\le x=\cos(\theta)\le 1 $, $\rho(z)\to e^{i\theta}\ $ (with
  $0\le\theta\le\pi)$,
so that $2a_n-1$ and $b_n$ behave like the {\it Fourier coefficients\/} of
order $2n$ and $2n+1$ of
a function whose singularities on $0<\theta<\pi$ are related to the
singularities of $w(\cos\theta)\sin\theta$. Singularities of type
$|\theta-\theta_0|^\gamma$ correspond indeed to $n^{-1-\gamma}
\cos(n\theta_0+$ const.) behaviour in the $n^{th}$ Fourier coefficient.

This analysis is no more valid for stronger singularity $w(x_0)=0$ or
$\infty$ because the approximations done near $x_0$ are no more
valid,
    subtle important effects take place in neighbourhoods of length about
$1/n$ of the singular points, see the famous [NevIII].

It is therefore not useless to have a close look at the simplest
orthogonal polynomials related to weights with interior singularities.
Special singular positions have been worked (in sieved polynomials
theory etc.), but here is something related to an \ \ {\it arbitrary\/}\quad
position:

\bigskip

{\bf 2.Freud's equations for the simplest generalized Jacobi polynomials.}

\bigskip

Let
$$\eqalign{
 w(x) &= B(1-x)^\alpha(x_0-x)^\gamma(1+x)^\beta\; x\in [-1,x_0]\;, \cr
  &=A(1-x)^\alpha(x-x_0)^\gamma(1+x)^\beta\; x\in [x_0,1]  \;, \cr
         }\eqno(3)
$$
with $-1<x_0<1$,  $A$ and $B >0$, $\alpha$, $\beta$ and $\gamma>-1$.

The case $A=B, \alpha=\beta, x_0=0$ has the simple solution deduced
from Jacobi polynomials
$a_n^2 =(n+2\alpha+\gamma\odd(n))(n+\gamma\odd(n))/
            [(2n+2\alpha+\gamma+1)
             (2n+2\alpha+\gamma-1)], b_n=0$, where
$\odd(n)=(1-(-1)^n)/2$. When $n\to\infty$, this case shows the asymptotic
behaviour $a_n\sim \dfrac12 - \dfrac{(-1)^n\gamma}{4n}, b_n=0$.

The $O(1/n)$ term is definitely related to the $|x|^\gamma$ behaviour of
the weight near 0, as shown by P.~Nevai (Theorem 4 of Section 7 of [Nev]):
if $w$ is even on $[-1,1]$, with $w(x)|x|^{-\gamma}$ positive and
continuously derivable on $(-1,1)$, then $a_n=1/2-\gamma (-1)^n/(4n)
+o(1/n)$ when $n\to\infty$.

We now try to investigate the recurrence coefficients $a_n$, $b_n$
when the weight is (3). This
weight
     is a {\it semi-classical\/} weight, as $w'/w$ is the same rational
function almost everywhere on the support $[-1,1]$ of $w$.

Semi-classical orthogonal polynomials have a rich differential structure,
according to a theory going as far as Laguerre [BeR, GaN, Lag, Mag2, Mag3,
Sho].

Freud [Fr] showed how to deduce recurrence coefficients asymptotics
from special identities. For a  general semi-classical weight
satisfying $w'(x)/w(x)=2V(x)/W(x)$ with $W(x)w(x)\to 0$ when $x$ tends to
any endpoint of the support $S$ of $w$ (Shohat's conditions [Sho]),
we find these identities
({\it Freud's equations\/}) by expanding
$$\eqalign{
 0 &= \int_S \left[W(x)w(x)p_n(x)p_{n-k}(x)\right]'\,dx \cr
   &= \int_S W(x)w(x)p'_n(x)p_{n-k}(x)\,dx
     +\int_S W(x)w(x)p_n(x)p'_{n-k}(x)\,dx +\cr
  &\hskip 50pt  +\int_S W'(x)w(x)p_n(x)p_{n-k}(x)\,dx
     +\int_S 2V(x)w(x)p_n(x)p_{n-k}(x)\,dx  \cr
  }$$
for $k=0,1$, remarking that any integral $\int_S P(x)w(x)p_n(x)p_{n-k}(x)dx$
where $P$ is a polynomial, is an expression involving $a_n, b_n,
a_{n\pm 1}, b_{n\pm 1},$ etc. according to $k$ and the degree of $P$,
and using $p'_n = n p_{n-1}/a_n +(b_0+\cdots +b_{n-1}-nb_{n-1})p_{n-2}/
(a_{n-1}a_n)+\cdots$ (see [BeR]).

Reduction to even measure:
one
considers the orthonormal polynomials $\{\tilde p_n\}$ with respect to the
even weight
$$ \tilde w(x) = 2|x|w(2x^2-1),\qquad {\rm for\ }
                     -1<x<1.$$
Then, $\tilde p_{2n}(x) =
p_n(2x^2-1)$, and one recovers the
recurrence relation for the $p_n$'s by contracting the recurrence
relation for the $\tilde p_n$'s:
$$\displaylines{
   \tilde a_{n+1}\tilde p_{n+1}(x) =x\tilde p_n(x) -\tilde a_n\tilde
    p_{n-1}(x) \qquad \Rightarrow  \hfill \cr
 \hfill \tilde a_{2n+1}\tilde a_{2n+2}\tilde p_{2n+2}(x) =
  (x^2 -\tilde a_{2n}^2-\tilde a_{2n+1}^2)\tilde p_{2n}(x)
     -  \tilde a_{2n-1}\tilde a_{2n}\tilde p_{2n-2}(x)\;,\cr
               }
$$
so:
$$a_n=2\tilde a_{2n-1} \tilde a_{2n}\;,\qquad
  b_n = -1+2\tilde a_{2n}^2+2\tilde a_{2n+1}^2 \eqno(4)$$
This allows to work with the single sequence $\{\tilde a_n\}$ instead
of the two sequences $\{a_n\}, \{b_n\}$.

Here,
$$\eqalign{
 \tilde w(x)=2|x|w(2x^2-1) &=
          \tilde B |x|^{2\beta+1}(\tilde x_0^2-x^2)^\gamma
                (1-x^2)^\alpha \qquad {\rm for\ } |x|<|\tilde x_0|,\cr
     &=   \tilde A |x|^{2\beta+1}(x^2-\tilde x_0^2)^\gamma
                (1-x^2)^\alpha \qquad {\rm for\ } |\tilde x_0|<|x|<1,\cr
          }
$$
where $\tilde x_0$ is the positive root of $2\tilde x_0^2-1=x_0$,
$\tilde A=2^{\alpha+\beta+\gamma+1}A$,
$\tilde B=2^{\alpha+\beta+\gamma+1}B$.

So, $\tilde W(x)=x(x^2-\tilde x_0^2)(x^2-1)=x^5-(\tilde x_0^2+1)x^3
+\tilde x_0^2 x$ and \hfill\break $2\tilde V(x) = (2\alpha+2\beta+2\gamma+1)x^4
-[2\alpha\tilde x_0^2 +(2\beta+1)(\tilde x_0^2+1)+2\gamma]x^2
+(2\beta+1)\tilde x_0^2$.

The equations for the $\tilde a_n$'s now follow from
Freud's method for
even weights, expanding $\tilde W\tilde w'=2\tilde V\tilde w$ as
$$\displaylines{
  \int_S \dfrac{\tilde W(x)}x \tilde w(x) \tilde p_n'(x)\tilde p_{n-1}(x)
 +\int_S \dfrac{\tilde W(x)}x \tilde w(x) \tilde p_n(x)\tilde p_{n-1}'(x)
 +\int_S \left(\dfrac{\tilde W(x)}x\right)'
                      \tilde w(x) \tilde p_n(x)\tilde p_{n-1}(x)+ \cr
  \hfill
 +\int_S \dfrac{2\tilde V(x)}x \tilde w(x) \tilde p_n(x)\tilde p_{n-1}(x)
=0,\cr}$$
using $\int_S \tilde w(x)\tilde p_i(x)\tilde p_j(x)dx=\delta_{i,j}$,
the recurrence relations (1) giving
      $\int_S x\tilde w(x)\tilde p_n(x)\tilde p_{n-1}(x)dx=\tilde a_n$,
  $\int_S x^2\tilde w(x)(\tilde p_n(x))^2 dx=\tilde a_n^2 +
                                             \tilde a_{n+1}^2$,
$\int_S \tilde w(x)x^{-1}\tilde p_n(x)\tilde p_{n-1}(x)dx=
\odd(n)/\tilde a_n$, etc., and
$\tilde p_n' = \dfrac{n}{\tilde a_n}\tilde p_{n-1} +
 \dfrac{2\sum_1^{n-1}\tilde a_k^2 -n\tilde a_{n-1}^2}
       {\tilde a_{n-2}\tilde a_{n-1}\tilde a_n} \tilde p_{n-3} +  $
\hfill\break
$ +
 \dfrac{n\tilde a_{n-3}^2 \tilde a_{n-1}^2
       -2(\tilde a_{n-3}^2+\tilde a_{n-2}^2+\tilde a_{n-1}^2)
         \sum_1^{n-1}\tilde a_k^2 +2\sum_1^{n-1}(\tilde a_k^4 +
                                          2\tilde a_k^2\tilde a_{k-1}^2)
       }{\tilde a_{n-4}
         \tilde a_{n-3}
         \tilde a_{n-2}
         \tilde a_{n-1} \tilde a_n } \tilde p_{n-5} +\cdots$\ ,
\hfill\break
one finally finds
$$\displaylines{
  2(n+\alpha+\beta+\gamma+2)\tilde a_n^2(\tilde a_{n-1}^2 +
                      \tilde a_n^2 + \tilde a_{n+1}^2)
     -2[\alpha\tilde x_0^2+(n+\beta+1)(\tilde x_0^2+1)+\gamma ]
           \tilde a_n^2   + \hfill\cr
    +2(2\tilde a_n^2-\tilde x_0^2-1)\sum_{j=1}^{n-1}\tilde a_j^2
   +n\tilde x_0^2 -2\tilde a_n^2\tilde a_{n-1}^2
   +2\sum_{j=1}^{n-1}(\tilde a_j^4+2\tilde a_j^2\tilde a_{j-1}^2)
   +(2\beta+1)\tilde x_0^2\,\odd(n) =0,  \cr
   \hfill n=1,2,\ldots \qquad   (5) \cr
               }
$$
($\tilde a_0=0$).

We see how any $\tilde a_n$ can be computed from the value of
$\tilde a_1$, which is the degree of freedom reflecting that
the same equations (5) hold for any choice of $\tilde A$ and
$\tilde B$ in the weight $\tilde w$. Actually $\tilde a_1$ is
linked to the ratio $\tilde A/\tilde B$ by
$$\eqalign{
   \tilde a_1^2 &= {\tilde\mu_2\over\tilde\mu_0} \cr &=
  {\tilde B \int_{|x|<\tilde x_0}|x|^{2\beta+3}(\tilde x_0^2-x^2)^\gamma
                (1-x^2)^\alpha\,dx +
   \tilde A \int_{|x|>\tilde x_0}|x|^{2\beta+3}(x^2-\tilde x_0^2)^\gamma
                (1-x^2)^\alpha\,dx  \over
   \tilde B \int_{|x|<\tilde x_0}|x|^{2\beta+1}(\tilde x_0^2-x^2)^\gamma
                (1-x^2)^\alpha\,dx +
   \tilde A \int_{|x|>\tilde x_0}|x|^{2\beta+1}(x^2-\tilde x_0^2)^\gamma
                (1-x^2)^\alpha\,dx
   }
        \cr }   \eqno(6)
$$

\medskip

%I wrote this on earlier drafts:

%``Asymptotic determination will be enormously simplified thanks to
%Szeg\H o's theory: the asymptotic behaviour of the apparently
%nasty sums is actually known up to the $O(1)$ term. Considering that
%the other terms of (3.1) contain $n$ factors, $O(1/n)$ accuracy will
%result!!!'', but things are not so simple after all.

Numerical experiments show that $\tilde a_2$, $\tilde a_3$,\dots can
be computed in a stable way from $\tilde a_1$ simply by considering
(5) as an equation for $\tilde a_{n+1}$ when $\tilde a_1$, \dots,
$\tilde a_n$ are known:
$$\displaylines{
   \tilde a_{n+1}^2 =
  \dfrac{\alpha\tilde x_0^2 +(n+\beta+1)(\tilde x_0^2+1)+\gamma}N
  -2\dfrac{\sum_1^{n-1}\tilde a_k^2}N +\cr
 +\dfrac{2(\tilde x_0^2+1)\sum_1^{n-1}\tilde a_k^2 -n\tilde x_0^2
 -2\sum_1^{n-1}(\tilde a_k^4+2\tilde a_k^2\tilde a_{k-1}^2)
 -(2\beta+1)\tilde x_0^2\odd(n)
         }{2N\tilde a_n^2}
  +\dfrac{\tilde a_{n-1}^2}N -\tilde a_n^2 -\tilde a_{n-1}^2  \cr
  \hfill(7)\cr
 }
$$
with $N=n+\alpha+\beta+\gamma+2$, $n=1,2,\ldots$.

\bigskip

{\bf 3.Asymptotic estimates.}

\bigskip

\rightline{Some people call these tricks ``special refinements'';}
\rightline{others call them  ``kludges''.\hskip 40pt}

\rightline{D.E.\ Knuth\ \hskip 75pt}
\bigskip

Putting an almost constant $\tilde a_n^2 \approx \tilde a^2$ in (5)
gives two possible asymptotic matches $\tilde a^2=1/4$ and
$\tilde a^2=\tilde x_0^2/4$, corresponding to weights with support
$[-1,1]$ and $[-\tilde x_0,\tilde x_0]$ (the latter when $A=0$).
More complicated behaviours are expected to hold when $\tilde x_0$ is
complex [GaN] and one should be able to establish correct asymptotic
behaviours from (5) alone, but this has not yet been achieved (the answer
to [GaN] was given in [N] with other techniques). Even when one knows
that $\tilde a_n^2\to 1/4$, some amount of guesswork will still be needed.
Let $\tilde a_n^2=\dfrac14+y_n$. We know from Szeg\H o's theory
([Sz] chap.\ 12)
that $\sum_1^{n-1}\tilde a_k^2=\dfrac{n}4+\xi+z_n$ and
     $\sum_1^{n-1}(\tilde a_k^4+2\tilde a_k^2\tilde a_{k-1}^2)=
       \dfrac{3n}{16}+\eta+u_n$ with $z_n$ and
 $u_n\to 0$ (See also [Nev2] p.91). So,
$$\displaylines{
   y_{n+1}=-\dfrac14+
  \dfrac{\alpha\tilde x_0^2 +(n+\beta+1)(\tilde x_0^2+1)+\gamma}N
  -2\dfrac{\dfrac{n}4+\xi+z_n}N  +\hfill\cr
 +\dfrac{2(\tilde x_0^2+1)(\dfrac{n}4+\xi+z_n)-n\tilde x_0^2
 -2(\dfrac{3n}{16}+\eta+u_n)
 -(2\beta+1)\tilde x_0^2\odd(n)
         }{2N(\dfrac14+y_n)}
  +\dfrac{\dfrac14+y_{n-1}}N -\dfrac12-y_n -y_{n-1}.\cr
  }
$$

We compute now $\xi$ and $\eta$:
from the Szeg\H o theory, let $\tilde\phi_n(z)=\tilde\kappa_n z^n +
\tilde\kappa'_n z^{n-2}+
\tilde\kappa'{}'_n z^{n-4}+\cdots$ be the orthonormal polynomials on the unit
circle with respect to $|\sin\theta|\tilde w(cos\theta) =
\tilde C(\theta)|\cos\theta|^{2\beta+1}
|\cos^2\theta-\cos^2(\theta_0/2)|^\gamma |\sin\theta|^{2\alpha+1}$,
with $\tilde C(\theta)=\tilde A$ on $|\theta|<  \theta_0/2$ and
                            $|\theta-\pi|<  \theta_0/2$ and
$\tilde C(\theta)=\tilde B$ elsewhere. The Szeg\H o function $\tilde D(z)$ whose
boundary values must be $|\tilde D(e^{i\theta})|=\sqrt{\tilde w(\cos\theta)
|\sin\theta|}$ is found by inspection to be
$$ \tilde D(z)=\tilde\kappa^{-1}(1-z^2)^{\alpha+1/2} (1+z^2)^{\beta+1/2}
             (1-e^{-i\theta_0}z^2)^{\gamma/2+i\lambda}
             (1-e^{i\theta_0}z^2)^{\gamma/2-i\lambda} ,$$
with $\tilde\kappa=2^{\alpha+\beta+\gamma+1}\tilde B^{(\theta_0-\pi)/(2\pi)}
                                            \tilde A^{-\theta_0/(2\pi)}$
and   $\lambda=(2\pi)^{-1}\log(\tilde B/\tilde A)$.
We find the limit values
of $\tilde\kappa'_n$ and $\tilde\kappa'{}'_n$ from the expansion of
 $1/\tilde D$:
$$\displaylines {
   \tilde\kappa'_n/\tilde\kappa_n \to \alpha-\beta+\gamma x_0
 +2\lambda\sin\theta_0,\cr
   \tilde\kappa'{}'_n/\tilde\kappa_n \to (\alpha-\beta)^2/2 +
   (\alpha+\beta+1)/2
 +(\alpha-\beta)(\gamma x_0+2\lambda\sin\theta_0)+\hfill\cr
 \hfill +(\gamma(\gamma+2)/4 -\lambda^2)\cos(2\theta_0)+
 \lambda(\gamma+1)\sin(2\theta_0)+\gamma^2/4+\lambda^2,\cr}$$
used in
$$\eqalign{
  \widetilde P_n(x)&=\dfrac{\tilde p_n(x)}{\tilde\gamma_n} =  \cr
    &= x^n -\left(\sum_1^{n-1}a_k^2\right)x^{n-2}
    +\left[\left(\sum_{k=1}^{n-1}\tilde a_k^2\right)^2 -
           \left(\sum_{k=1}^{n-1}\tilde a_k^4\right) -2
           \left(\sum_{k=1}^{n-2}\tilde a_k^2\tilde a_{k+1}^2\right)\right]
        x^{n-4}/2+\cdots \cr
 &= \dfrac{z^{-n}\tilde\phi_{2n}(z)+z^n\tilde\phi_{2n}(z^{-1})}
               {2^n(\tilde\kappa_{2n}+\tilde\phi_{2n}(0)} \cr
  &\sim \dfrac{T_n(x)}{2^{n-1}}+
    \dfrac{\tilde\kappa'}{\tilde\kappa}  \dfrac{T_{n-2}(x)}{2^{n-1}}+
    \dfrac{\tilde\kappa'{}'}{\tilde\kappa}  \dfrac{T_{n-4}(x)}{2^{n-1}}+
    \cdots \cr
  &\sim x^n -\dfrac{n-\tilde\kappa'/\tilde\kappa}4\, x^{n-2}
+\dfrac{n(n-3)/2 -(n-2)\tilde\kappa'/\kappa +\tilde\kappa'{}'/\tilde\kappa}
       {16} x^{n-4}+\cdots
  \cr  }$$
($z+z^{-1}=2x$).

So,
$$ \xi = \lim_{n\to\infty} \sum_1^{n-1}\tilde a_k^2 -n/4 =
      -\kappa'/(4\kappa) =
      -(\alpha-\beta+\gamma x_0 +2\lambda\sin\theta_0)/4, \eqno(8) $$
$\eta=\lim_{n\to\infty}\sum_1^{n-1}
                     (\tilde a_k^4+2\tilde a_k^2\tilde a_{k-1}^2)
 -3n/16 = ((\kappa'/\kappa)^2 -4\kappa'/\kappa-2\kappa'{}'/\kappa)/16$.

The equation for the $y_n$'s reduces to
$$\displaylines{
  y_{n+1}-2x_0 y_n+y_{n-1} =
  \dfrac{(x_0+1)(\beta+1/2)(-1)^n +2(x_0+2)z_n +y_{n-1}-4u_n}N  -\cr
  -\dfrac{16y_n}{1+4y_n}
   \dfrac{(2x_0+1)/4 -\lambda\sin\theta_0 +(x_0+1)(\beta+1/2)(-1)^n
    +2(x_0+3)z_n-4u_n}{4N}\cr
     \hfill -\dfrac{4(2x_0+1)y_n^2}{1+4y_n}\qquad(9)\cr} $$
where $z_n=-\sum_n^\infty y_k$ and $u_n=-y_{n-1}/2-\sum_n^\infty(
   3y_k/2+y_k^2+2y_ky_{k-1})$.

%This is the big simplification I referred to above.
The form (9) should give hints on the behaviour of $y_n$ when $n\to\infty$.
No proof will be attempted here, only reasonable asymptotic matching
and numerical checks. Use of Painlev\'e-like differential equations in
$\theta_0$ is another method of investigation which could be used in the
future (see [Mag3]).

 The right-hand side is
small, even with respect to the $y$'s, so $y_{n+1}-2x_0y_n+y_{n-1}$
is small, and this suggests a $\exp(\pm in\theta_0)$ behaviour somewhere.
However, I still don't have a tight proof that $y_n$, $z_n$ and $u_n$
are $O(1/n)$. Assuming $y_n=K_1 (-1)^n /n +K_2 e^{in\theta_0}/n^{\zeta_2}
                                          +K_3 e^{-in\theta_0}/n^{\zeta_3}$,
matching the two sides gives $K_1=-\beta/2-1/4$, $\zeta_{2,3}=1\pm 2i\lambda$.
Numerical checks have been performed on the form
$$\tilde a_n^2-1/4=y_n=    -(\beta+1/2)(-1)^n/(2n)
  +K \cos(n\theta_0 -2\lambda\log n -\varphi)/n +o(1/n)\eqno(10)$$
where $K$ and $\varphi$ are unknown functions of $\alpha$, $\beta$,
$\gamma$, $\lambda$ and $x_0$. Given $\alpha$, $\beta$, $\gamma$ and
$x_0$, the algorithm first performs (7) with several trial starting
values $\tilde a_1^2$ and computes the corresponding $\lambda$ from (8),
allowing the determination of the coefficients in (6). It is then
possible to run (7) for a requested value of $\lambda$, and to estimate
$K$ and $\varphi$ in (10) from numerical values of $\tilde a_n^2$ for
large $n$ (up to the 10000-100000 range). Very satisfactory empirical
formulas for $K$ and $\varphi$ appear to be
 $K= (\gamma^2/4+\lambda^2)^{1/2} \sin(\theta_0/2)$ and
$\varphi=(\alpha+1+\gamma/2)\pi-
     (\alpha+\beta+\gamma+1/2)\theta_0+2\lambda \log(2\sin \theta_0)
     -2\arg\Gamma(\gamma/2+i\lambda) -\arg(\gamma/2+i\lambda)$.

Whence, from (4): $a_n-1/2\sim y_{2n-1}+y_{2n}, b_n\sim
2(y_{2n}+y_{2n+1})$,  the

\bigskip

{\bf Conjecture.\/} The recurrence coefficients related to the
simplest generalized Jacobi weight (3) satisfy

$$ \eqalign{
    a_n &= \dfrac12 -\dfrac{M}n
                \cos\left[2n\theta_0-2\lambda\log(4n\sin\theta_0)
                     -\Phi\right]
         +o\left(\dfrac1{n}\right) ,\cr
  b_n &=  -\dfrac{2M}n \cos\left[(2n+1)\theta_0-2\lambda\log(4n\sin\theta_0)
                     -\Phi\right]
         +o\left(\dfrac1{n}\right) ,\cr
           }
$$
when $n\to\infty$, where $x_0=\cos\theta_0$, $0<\theta_0<\pi$,
   $\lambda=\log(B/A)/(2\pi)$,
 $M=\dfrac12 (\gamma^2/4+\lambda^2)^{1/2}\sin\theta_0$,
$\Phi=(\alpha+\gamma/2)\pi-(\alpha+\beta+\gamma)\theta_0
     -2\arg\Gamma(\gamma/2+i\lambda) -\arg(\gamma/2+i\lambda)$.

\bigskip

Here is a sample of the numerical check: $K$ and $\varphi$ are
extracted from the form (10) on a sample of $\tilde a_n^2$, with $n$
going up to 500000. The values of $K$ and $\lambda$ (from (8)) are quite
stable, but things are not so easy with $\varphi$ (the ``phase'' column).
Finally, the line ``check'' contains the computed values of $K$
and $\varphi$, and a new check of $\lambda$ through Tur\'an
determinant weight function reconstruction yielding $\tilde A$ and
$\tilde B$.

\bigskip

\verblist{magnus1.l}

\bigskip

{\bf Acknowledgements.}

Many thanks to H.Dette, L.Golinskii and P.Nevai for their interest.
The computations have been performed on the Convex {\tt C3820\/} of
the University.

% if this happens to be accepted for the Gatteschi proceedings:
Many thanks to G.~Allasia who manages to edit the present proceedings.

\bigskip
{\bf References\/}

\frenchspacing      \parindent=40pt

\item{Ap} A.I.\ APTEKAREV, Asymptotics of orthogonal polynomials in a
     neighborhood of the endpoints of the interval of orthogonality,
    {\sl Russ.\ Acad.\ Nauk Mat.\ Sb.\/} {\bf 183\/} (1992) =
    {\sl Russian Acad.\ Sci.\ Sb.\ Math.\/} {\bf 76\/} (1993)  35-50.

\item{Ba} W.\ BARRETT, An asymptotic formula relating to orthogonal
         polynomials, {\sl J.\ London Math.\ Soc.\ (2)\/}
         {\bf 6\/} (1973), 701-704.

\item{BeR}  S. BELMEHDI, A. RONVEAUX, About non linear systems
     satisfied by the recurrence coefficients of semiclassical
     orthogonal polynomials, to appear in {\sl J.\ Approx.\ Th.\/}

\item{Fr}  G.FREUD, On the coefficients in the recursion formul\ae\  of
     orthogonal polynomials, {\sl Proc.\ Royal Irish Acad. Sect. A\/}
     {\bf 76\/} (1976), 1-6.

\item{GaN} J.L.~GAMMEL, J.~NUTTALL, Note on generalized Jacobi
     polynomials, {\it in\/} ``The Riemann Problem, Complete
     Integrability and Arithmetic Applications'' (D.~Chudnovsky and
     G.~Chudnovski, Eds.), pp.258-270, Springer-Verlag (Lecture
     Notes Math. {\bf 925\/}), Berlin, 1982.

\item{Ho} C.H.~HODGES, Van Hove singularities and continued fraction
     coefficients, {\sl J.\ Physique Lett.\ \/} {\bf 38\/} (1977),
    L187-L189.

\item{Lag}  E. LAGUERRE, Sur la r\'eduction en fractions continues d'une
      fraction qui satisfait \`a une \'equation diff\'erentielle lin\'eaire
      du premier ordre dont les coefficients sont rationnels,
      {\sl J. Math. Pures Appl. (4)\/} {\bf 1\/} (1885), 135-165 =
      pp. 685-711 {\it in\/}
      {\sl Oeuvres\/}, Vol.II, Chelsea, New-York 1972.

\item{LaG} Ph.\ LAMBIN, J.P.\ GASPARD, Continued-fraction technique for
      tight-binding systems. A generalized-moments approach,
      {\sl Phys.\ Rev.\ B\/} {\bf 26\/} (1982) 4356-4368.

\item  {Mag1}  A.P.~MAGNUS, Recurrence coefficients for orthogonal
     polynomials on connected and non connected sets, {\it in\/}
     ``Pad\'e Approximation and its Applications, Proceedings,
     Antwerp 1979'' (L.~Wuytack, editor), pp.~150-171, Springer-Verlag
     (Lecture Notes Math.\ {\bf 765\/}), Berlin, 1979.

\item{Mag2} A.P.~MAGNUS, On Freud's equations for exponential
       weights, {\sl J.\ Approx.\ Th.\/} {\bf 46\/} (1986) 65-99.

\item{Mag3} A.P.~MAGNUS, Painlev\'e-type differential equations for
        the recurrence  coefficients of semi-classical orthogonal
         polynomials, to appear in {\sl J.\ Comp.\ Appl.\ Math.}

\item{Mart} A.MARTIN {\it et al.\/}, L\'eon Van Hove, 1924-1990,
      {\sl Courrier CERN\/} {\bf 30\/} n$\rondsur{\ } 2$ (March 1991),
       20-27.

\item{Nev} P.NEVAI, {\sl Orthogonal Polynomials, Memoirs AMS\/}
      vol.\ {\bf 18\/} nr.\ {\bf 213\/} (March 1979).

\item{NeIII} P.NEVAI, Mean convergence of Lagrange interpolation III,
     {\sl Trans.\ Amer.\ Math.\ Soc.\/} {\bf 282\/} (1984) 669-698.

\item{Nev2} P.\ NEVAI, Orthogonal polynomials, recurrences, Jacobi
     matrices and measures, pp. 79-104 {\it in\/}
     {\sl Progress in Approximation Theory, an International
     Perspective\/} (A.A.GONCHAR \& E.B.SAFF, editors),
     Springer-Verlag, New York, 1992.

\item{NV} P.NEVAI, W.~VAN~ASSCHE, Compact perturbations of orthogonal
      polynomials, {\sl Pacific J.\ Math.\/} {\bf 153\/} (1992) 163-184.

\item {N} J.NUTTALL, Asymptotics of generalized Jacobi polynomials,
      {\sl Constr.\ Approx.\/} {\bf 2\/} (1986) 59-77.

\item {Sho}  J.A. SHOHAT, A differential equation for orthogonal
     polynomials, {\sl Duke Math. J.\/} {\bf 5\/} (1939),401-417.

\item {Sz}  G. SZEG\H O, {\sl Orthogonal Polynomials\/},
          Colloquium Publications, Vol. 23
      Amer. Math. Soc.  Providence, Rhode Island,   1967.

\item {VA} W.~VAN~ASSCHE, Lecture Notes, Namur, Feb.\ 1993.

\bye